\documentclass[a4paper,12pt,reqno]{amsart}
\usepackage{amssymb}
\usepackage{amsmath}
\usepackage{ifthen}
\usepackage[dvips]{graphicx}
\nonstopmode \numberwithin{equation}{section}
\usepackage{amssymb}
\usepackage{ifthen}
\usepackage{graphicx}
\usepackage{amsmath}
\usepackage[T1]{fontenc} 
\usepackage{comment}

\nonstopmode \numberwithin{equation}{section}
\theoremstyle{plain}
\newtheorem{thm}{Theorem}
\numberwithin{thm}{section}
\newtheorem{cor}{Corollary}
\numberwithin{cor}{section}
\newtheorem{lem}{Lemma}
\numberwithin{lem}{section}
\newtheorem{prop}{Proposition}

\newtheorem{conj}{Conjecture}

\newcommand\numberthis{\addtocounter{equation}{1}\tag{\theequation}}
\theoremstyle{definition}
\newtheorem{defn}{Definition}[section]

\newtheorem{prob}{Problem}
\newtheorem{rem}{Remark}[section]

\newcounter{minutes}
\divide\time by 60
\newcounter{hours}
\multiply\time by 60
\addtocounter{minutes}{-\time}

\newcounter {own}
\def\theown {\thesection       .\arabic{own}}

\newenvironment{pf}[1][]{%
	\vskip 3mm
	\noindent
	\ifthenelse{\equal{#1}{}}%
	{{\slshape Proof. }}%
	{{\slshape #1.} }%
}%
{\qed\bigskip}

\theoremstyle{plain}
\newtheorem{Thm}{Theorem}

\newcommand{\ID}{{\mathbb D}}

\newcommand{\D}{{\mathbb D}}

\renewcommand{\theequation}{\thesection.
\arabic{equation}}
\numberwithin{equation}{section}
\def\be{\begin{equation}}
\def\ee{\end{equation}}

\newcommand{\bee}{\begin{enumerate}}
	\newcommand{\eee}{\end{enumerate}}

\newcommand{\blem}{\begin{lem}}
	\newcommand{\elem}{\end{lem}}
\newcommand{\bthm}{\begin{thm}}
	\newcommand{\ethm}{\end{thm}}
\newcommand{\bcor}{\begin{cor}}
	\newcommand{\ecor}{\end{cor}}
\newcommand{\beg}{\begin{examp}}
	\newcommand{\eeg}{\end{examp}}
\newcommand{\begs}{\begin{examples}}
	\newcommand{\eegs}{\end{examples}}
\allowdisplaybreaks
\newcommand{\bdefn}{\begin{defn}}
	\newcommand{\edefn}{\end{defn}}

\newcommand{\bprob}{\begin{prob}}
	\newcommand{\eprob}{\end{prob}}
\newcommand{\bei}{\begin{itemize}}
	\newcommand{\eei}{\end{itemize}}

\newcommand{\bcon}{\begin{conj}}
	\newcommand{\econ}{\end{conj}}
\newcommand{\bcons}{\begin{conjs}}
	\newcommand{\econs}{\end{conjs}}
\newcommand{\bprop}{\begin{prop}}
	\newcommand{\eprop}{\end{prop}}
\newcommand{\br}{\begin{rem}}
	\newcommand{\er}{\end{rem}}
\newcommand{\brs}{\begin{rems}}
	\newcommand{\ers}{\end{rems}}
\newcommand{\bo}{\begin{obser}}
	\newcommand{\eo}{\end{obser}}
\newcommand{\bos}{\begin{obsers}}
	\newcommand{\eos}{\end{obsers}}
\newcommand{\bpf}{\begin{pf}}
	\newcommand{\epf}{\end{pf}}
\newcommand{\ba}{\begin{array}}
	\newcommand{\ea}{\end{array}}
\newcommand{\beq}{\begin{eqnarray}}
\newcommand{\beqq}{\begin{eqnarray*}}
\newcommand{\eeq}{\end{eqnarray}}
\newcommand{\eeqq}{\end{eqnarray*}}

\begin{document}

\title{Landau type theorem for $\alpha$-harmonic mappings}

\author{Vasudevarao Allu}
\address{Vasudevarao Allu,
	School of Basic Science,
	Indian Institute of Technology Bhubaneswar,
	Bhubaneswar-752050, Odisha, India.}
\email{avrao@iitbbs.ac.in}

\author{Rohit Kumar}
\address{Rohit Kumar,
	School of Basic Science,
	Indian Institute of Technology Bhubaneswar,
	Bhubaneswar-752050, Odisha, India.}
\email{rohitk12798@gmail.com}

\subjclass[{AMS} Subject Classification:]{Primary 31A30, 30C99; Secondary 31A05, 31A35.}

\def\thefootnote{}
\footnotetext{ {\tiny File:~\jobname.tex,
		printed: \number\year-\number\month-\number\day,
		\thehours.\ifnum\theminutes<10{0}\fi\theminutes }
} \makeatletter\def\thefootnote{\@arabic\c@footnote}\makeatother

\begin{abstract}
In this paper, we first obtain an estimate of the coefficients for $\alpha$-harmonic mappings. By applying these coefficient estimates, we prove the Landau type theorem for $\alpha$-harmonic mappings defined on the unit disc $\ID$.

\end{abstract}

\maketitle
\pagestyle{myheadings}
\markboth{Vasudevarao Allu and  Rohit Kumar}{Landau type theorem for $\alpha$-harmonic mappings}

\section{\textbf{Introduction}}
Let $\mathbb{D}=\{z\in\mathbb{C}:|z|<1\}$ represents the unit disc in the complex plane $\mathbb{C}$. We denote the complex differential operators 
$$
\frac{\partial}{\partial z}=\frac{1}{2}\left(\frac{\partial}{\partial x}-i \frac{\partial}{\partial y}\right) \text { and } \frac{\partial}{\partial \bar{z}}=\frac{1}{2}\left(\frac{\partial}{\partial x}+i \frac{\partial}{\partial y}\right)
$$
for $z=x+i y \in \mathbb{C}$ where $x$ and $y$ are real.
For a continuously differentiable function $f$, denote
$$
\Lambda_f=\max _{0 \leq \theta \leq 2 \pi}\left|f_z+e^{-2 i \theta} f_{\bar{z}}\right|=\left|f_z\right|+\left|f_{\bar{z}}\right|$$ and
$$
\lambda_f=\min _{0 \leq \theta \leq 2 \pi}\left|f_z+e^{-2 i \theta} f_{\bar{z}}\right|=|| f_z|-| f_{\bar{z}}||,
$$
where $f_z=\partial f / \partial z$ and $f_{\bar{z}}=\partial f / \partial \bar{z}$. Throughout this paper, we denote by $C^n(\mathbb{D})$ the set of all $n$-times continuously differentiable complex-valued functions in $\mathbb{D}$, where $n \in\{1,2, \cdots\}$.
\vspace{3mm}

For $\alpha \in \mathbb{R}$ and $z \in \ID$, let
$$
T_\alpha=-\frac{\alpha^2}{4}\left(1-|z|^2\right)^{-\alpha-1}+\frac{\alpha}{2}\left(1-|z|^2\right)^{-\alpha-1}\left(z \frac{\partial}{\partial z}+\bar{z} \frac{\partial}{\partial \bar{z}}\right)+\frac{1}{4}\left(1-|z|^2\right)^{-\alpha} \Delta
$$
be the second order elliptic partial differential operator, where $\Delta$ is the usual complex Laplacian operator
$$
\Delta:=4 \frac{\partial^2}{\partial z \partial \bar{z}}=\frac{\partial^2}{\partial x^2}+\frac{\partial^2}{\partial y^2} .
$$
The corresponding homogeneous differential equation is
\begin{equation}\label{Rohi-Vasu-P3-equation-001}
T_\alpha(f)=0 \quad \text { in } \mathbb{D} \text {. }
\end{equation}
and its associated Dirichlet boundary value problem is as follows
\begin{equation}\label{Rohi-Vasu-P3-equation-002}
\begin{cases}T_\alpha(f)=0 & \text { in } \mathbb{D}, \\ f=f^* & \text { on } \partial \mathbb{D}. \end{cases}
\end{equation}
Here, the boundary data $f^* \in \mathfrak{D}^{\prime}(\partial \mathbb{D})$ is a distribution on the boundary $\partial \mathbb{D}$ of $\mathbb{D}$, and the boundary condition in (\ref{Rohi-Vasu-P3-equation-002}) is interpreted in the distributional sense that $f_r \rightarrow f^*$ in $\mathfrak{D}^{\prime}(\partial \mathbb{D})$ as $r \rightarrow 1^-$, where
$$
f_r(e^{i \theta})=f\left(r e^{i \theta}\right), e^{i \theta} \in \partial \mathbb{D},
$$
for $r \in[0,1)$. In 2014, Olofsson \cite{Olofsson-2014} proved that, for parameter values $\alpha>-1$, if a function $f \in C^2(\mathbb{D})$ satisfies (\ref{Rohi-Vasu-P3-equation-001}) with $\lim _{r \rightarrow 1^-} f_r=f^* \in \mathfrak{D}^{\prime}(\partial \mathbb{D})$, then it has the form of a Poisson type integral
$$
f(z)=\frac{1}{2 \pi} \int_0^{2 \pi} K_\alpha\left(z e^{-i \tau}\right) f^*\left(e^{i \tau}\right) d\tau \quad \text {for } z \in \mathbb{D},
$$
where
$$
K_\alpha(z)=c_\alpha \frac{\left(1-|z|^2\right)^{\alpha+1}}{|1-z|^{\alpha+2}},
$$
$c_\alpha=(\Gamma(\alpha / 2+1))^2 / \Gamma(1+\alpha)$ and $\Gamma(s)=\int_0^{\infty} t^{s-1} e^{-t} d t$ for $s>0$ is the standard Gamma function. 
\vspace{2mm}

If we take $\alpha=2(n-1)$, then $f$ is polyharmonic (or $n$-harmonic), where $n \in\{1,2,3, \ldots\}$. For a detailed study on polyharmonic mappings, we refer to \cite{Borichev-Hedenmalm-2014,Chen-Vuorinen-2015,Li-Ponnusamy-2018,Li-Wang-Xiao-2017,Long-Wang-2021,Long-Wang-2022,Olofsson-2005, Olofsson-2014}. In particular, if $\alpha=0$, then $f$ is harmonic. Thus, $f$ is a kind of generalization of classical harmonic mappings.
\vspace{2mm}
 
The Gauss hypergeometric function is defined by the following series
$$
F(a, b ; c ; x)=\sum_{n=0}^{\infty} \frac{(a)_n(b)_n}{(c)_n} \frac{x^n}{n !}
$$
for $|x|<1$, where $(a)_0=1$ and $(a)_n=a(a+1) \cdots(a+n-1)$ for $n=1,2, \ldots$ are the Pochhammer symbols. Obviously, for $n=0,1,2, \ldots,(a)_n=\Gamma(a+n) / \Gamma(a)$. It is easy to see that
$$
\frac{d}{d x} F(a, b ; c ; x)=\frac{a b}{c} F(a+1, b+1 ; c+1 ; x) .
$$
Furthermore, for $\operatorname{Re}(c-a-b)>0$, we have (see \cite[Theorem 2.2.2]{Andrews-Askey-Roy-1999})
\begin{equation}\label{Rohi-Vasu-P3-equation-003}
F(a, b ; c ; 1)=\lim _{x \rightarrow 1} F(a, b ; c ; x)=\frac{\Gamma(c) \Gamma(c-a-b)}{\Gamma(c-a) \Gamma(c-b)}.
\end{equation}
 The following result concerns the solution to the equation (\ref{Rohi-Vasu-P3-equation-001}).
\begin{Thm}\cite{Olofsson-2014}\label{Thm-A}
 Let $\alpha \in \mathbb{R}$ and $f \in C^2(\mathbb{D})$. Then $f$ satisfies (\ref{Rohi-Vasu-P3-equation-001}) if, and only if, it has a series expansion of the form
\begin{equation}\label{Rohi-Vasu-P3-equation-003-a}
f(z)=\sum_{k=0}^{\infty} c_k F\left(-\frac{\alpha}{2}, k-\frac{\alpha}{2} ; k+1 ;|z|^2\right) z^k +\sum_{k=1}^{\infty}{c_{-k}} F\left(-\frac{\alpha}{2}, k-\frac{\alpha}{2} ; k+1 ;|z|^2\right) \bar{z}^k, 
\end{equation}
for some sequence $\left\{c_k\right\}_{-\infty}^{\infty}$ of complex number satisfying
\begin{equation}\label{Rohi-Vasu-P3-equation-003-b}
\lim _{|k| \rightarrow \infty} \sup \left|c_k\right|^{\frac{1}{|k|}} \leq 1 .
\end{equation}
In particular, the expansion (\ref{Rohi-Vasu-P3-equation-003-a}), subject to the condition (\ref{Rohi-Vasu-P3-equation-003-b}), converges in $C^{\infty}(\ID)$, and every solution $f$ of (\ref{Rohi-Vasu-P3-equation-001}) is $C^{\infty}(\ID)$ in the unit disc $\ID$.
\end{Thm}
In 2014, Olofsson \cite{Olofsson-2014} pointed out that if $\alpha \leq-1, f \in C^2(\mathbb{D})$ satisfies (1.1), and the boundary limit $f^*=$ $\lim _{r \rightarrow 1^{-}} f_r$ exists in $\mathfrak{D}^{\prime}(\partial \mathbb{D})$, then $f(z)=0$ for all $z \in \mathbb{D}$. Therefore, in this paper, we always assume that $\alpha>-1$.
\vspace{2mm}

The classical Landau theorem (see \cite{Landau-1926}) asserts that if $f$ is a holomorphic mapping with $f(0)=0$, $f'(0)=1$ and $|f(z)|<M$ for $z \in \mathbb{D}$, then $f$ is univalent in $\mathbb{D}_{r_0}$, and $f(\mathbb{D}_{r_0})$ contains a disc $\mathbb{D}_{\sigma_0}$, where $$r_0=\frac{1}{M+\sqrt{M^2-1}}\ \ \text{and} \ \ \sigma_0= Mr_0^2.$$
  The quantities $r_0$ and $\sigma_0$ can't be improved as the function $f_0(z)= Mz \left(\frac{1-Mz}{M-z}\right)$ shows the sharpness of the quantities $r_0$ and $\sigma_0$. 
\vspace{2mm}

In 2000, Chen {\it{et al.}} \cite{Chen-Gauthier-Hengertner-2000} obtained the Landau type theorems for harmonic mappings in the unit disc $\ID$. 

\begin{Thm}\cite{Chen-Gauthier-Hengertner-2000}\label{Thm-B}
Let $f$ be a harmonic mapping of the unit disc $\mathbb{D}$ such that $f(0)=0$, $f_{\bar{z}}(0)=0, f_z(0)=1$, and $|f(z)|<M$ for $z \in \mathbb{D}$. Then, $f$ is univalent on a disc $\mathbb{D}_{\rho_0}$ with
$$
\rho_0=\frac{\pi^2}{16 m M},
$$
and $f\left(\mathbb{D}_{\rho_0}\right)$ contains a schlicht disc $\mathbb{D}_{R_0}$ with
$$
R_0=\rho_0 / 2=\frac{\pi^2}{32 m M},
$$
where $m \approx 6.85$ is the minimum of the function $(3-r^2)/(r(1-r^2))$ for $0<r<1$.
\end{Thm}

\begin{Thm}\cite{Chen-Gauthier-Hengertner-2000}\label{Thm-C}
 Let $f$ be a harmonic mapping of the unit disc $\mathbb{D}$ such that $f(0)=0$, $\lambda_f(0)=1$ and $\Lambda_f(z) \leq \Lambda$ for $z \in \mathbb{D}$. Then, $f$ is univalent on a disc $\mathbb{D}_{\rho_0}$ with
$$
\rho_0=\frac{\pi}{4(1+\Lambda)},
$$
and $f\left(\mathbb{D}_{\rho_0}\right)$ contains a schlicht disc $\mathbb{D}_{R_0}$ with
$$
R_0=\frac{1}{2} \rho_0=\frac{\pi}{8(1+\Lambda)} .
$$
\end{Thm}
The radii $\rho_o$ and $R_0$ in Theorem \ref{Thm-B} and Theorem \ref{Thm-C} are not sharp. These results have been improved by Chen {\it{et al.}} \cite{Chen-Ponnusamy-Wang-2011}, Dorff and Nowak \cite{Dorff-Nowak-2004}, Grigoryan \cite{Grigoyan-2006}, Huang \cite{Huang-2008}, Liu \cite{M. Liu-2009, Liu-2009}, Liu and Chen \cite{Liu-Chen-2018} and Zhu \cite{Zhu-2015}. In particular, Liu \cite{Liu-2009} has proved that under the hypothesis of Theorem \ref{Thm-C}, $\Lambda\geq 1$, and when $\Lambda=1$, $f$ is univalent in the disc $\ID_{\rho}$ and $f(\ID_{\rho})$ contains a schlicht disc $\ID_R$ with $\rho = R = 1$ being sharp.  In 2018, Liu and Chen \cite{Liu-Chen-2018} established the sharp version of Theorem \ref{Thm-C} for $\Lambda>1$ by applying geometric method. The Landau type theorem has also been studied for elliptic harmonic mappings \cite{Allu-Kumar-2024}, logharmonic mappings \cite{Abdulhadi-Muhanna-Ali-2012} and log-$p$-harmonic mappings \cite{Li-Wang-2012}.
\vspace{2mm}

 In 2008, Abdulhadi and Muhanna \cite{Abdulhadi-Muhanna-2008} established two versions of Landau-type theorems for certain bounded biharmonic mappings. After that, the Landau type theorem for functions $f$ satisfying (\ref{Rohi-Vasu-P3-equation-001}) with $\alpha=2(n-1)$,  where $n \in \{1,2,3,\ldots\}$, has been extensively studied by many authors (see \cite{Abdulhadi-Muhanna-2008,Chen-Ponnusamy-Wang-2009,Chen-Rasila-Wang-2014,Liu-2008,Zhu-Liu-2013}).

\section{Coefficient Estimates}
There have been many discussions in the literature regarding the coefficient estimates of the functions $f$ satisfying (\ref{Rohi-Vasu-P3-equation-001})  with $\alpha =0$ and $\alpha=2.$ In 2015, Chen and Vuorinen \cite{Chen-Vuorinen-2015} proved the coefficient estimates for bounded functions satisfying (\ref{Rohi-Vasu-P3-equation-001}) with $\alpha>-1$. In this section, we prove a different coefficient estimates for $\alpha$-harmonic mappings when $\alpha>0$.
\begin{thm}\label{Rohi-Vasu-P3-Theorem-001}
For $\alpha>0$, let $f\in C^2(\ID)$ be the solution of (\ref{Rohi-Vasu-P3-equation-001}) with series expansion of the form (\ref{Rohi-Vasu-P3-equation-003-a}) and $\Lambda_f(z)\leq \Lambda$ for each $z\in \ID$, where $\Lambda$ is a positive number. Then, for $k\in \{1,2,3,...\}$,
\begin{align*}
\frac{\alpha|2k-\alpha|}{2(k+1)} &\biggl(\left|c_k F\left(-\frac{\alpha}{2}+1, k-\frac{\alpha}{2}+1 ; k+2 ;1\right)\right|\\
&+ \left|c_{-k} F\left(-\frac{\alpha}{2}+1, k-\frac{\alpha}{2}+1 ; k+2 ;1\right)\right|\biggl)\leq \Lambda.
\end{align*}

\end{thm}

\begin{pf}
We know that 
\begin{equation}\label{Rohi-Vasu-P3-equation-003-c}
F\left(-\frac{\alpha}{2}, k-\frac{\alpha}{2} ; k+1 ;|z|^2\right)=\sum_{n=0}^{\infty}\frac{(-\frac{\alpha}{2})_n(k-\frac{\alpha}{2})_n}{(k+1)_n}\frac{|z|^{2n}}{n!}.
\end{equation}
Differentiating (\ref{Rohi-Vasu-P3-equation-003-c}) partially with respect to $z$, we obtain
\begin{align*}
\frac{\partial}{\partial z}F\left(-\frac{\alpha}{2}, k-\frac{\alpha}{2} ; k+1 ;|z|^2\right)&=\bar{z}\sum_{n=1}^{\infty}\frac{(-\frac{\alpha}{2})_n(k-\frac{\alpha}{2})_n}{(k+1)_n}\frac{|z|^{2(n-1)}}{(n-1)!}\\
&= \bar{z}\frac{-\frac{\alpha}{2}( k-\frac{\alpha}{2})}{k+1}F\left(-\frac{\alpha}{2}+1, k-\frac{\alpha}{2}+1 ; k+2 ;|z|^2\right).
\end{align*}
Similarly, we have
$$
\frac{\partial}{\partial \bar{z}}F\left(-\frac{\alpha}{2}, k-\frac{\alpha}{2} ; k+1 ;|z|^2\right) = z\frac{-\frac{\alpha}{2}( k-\frac{\alpha}{2})}{k+1}F\left(-\frac{\alpha}{2}+1, k-\frac{\alpha}{2}+1 ; k+2 ;|z|^2\right).
$$
For $r\in [0,1)$, let 
$$
G_k(r,\alpha)=\frac{-\frac{\alpha}{2}( k-\frac{\alpha}{2})}{k+1}F\left(-\frac{\alpha}{2}+1, k-\frac{\alpha}{2}+1 ; k+2 ;r^2\right).
$$
A simple computation shows that 
\begin{align*}\label{Rohi-Vasu-P3-equation-004}
f_z(re^{i\theta})=\sum_{k=1}^{\infty}c_kr^{k+1}G_k(r,\alpha)e^{i(k-1)\theta}&+\sum_{k=1}^{\infty} k c_k F\left(-\frac{\alpha}{2}, k-\frac{\alpha}{2} ; k+1 ;r^2\right)r^{k-1}e^{i(k-1)\theta}\\ &+\sum_{k=1}^{\infty}c_{-k}r^{k+1}G_k(r,\alpha)e^{-i(k+1)\theta}\numberthis
\end{align*}
and 
\begin{align*}\label{Rohi-Vasu-P3-equation-005}
f_{\bar{z}}(re^{i\theta})=\sum_{k=1}^{\infty}c_kr^{k+1}G_k(r,\alpha)e^{i(k+1)\theta}&+\sum_{k=1}^{\infty} k c_{-k} F\left(-\frac{\alpha}{2}, k-\frac{\alpha}{2} ; k+1 ;r^2\right)r^{k-1}e^{-i(k-1)\theta}\\ 
&+\sum_{k=1}^{\infty}c_{-k}r^{k+1}G_k(r,\alpha)e^{-i(k-1)\theta},\numberthis
\end{align*}
where $z=re^{i\theta}$. By multiplying (\ref{Rohi-Vasu-P3-equation-004}) by $e^{i(k+1)\theta}$ and then integrating from $0$ to $2\pi$, we obtain
\begin{equation}\label{Rohi-Vasu-P3-equation-006}
c_{-k} r^{k+1}G_k(r,\alpha)=\frac{1}{2\pi}\int_0^{2\pi}f_z(re^{i\theta}) e^{i(k+1)\theta}d\theta.
\end{equation}
Similarly, multiplying (\ref{Rohi-Vasu-P3-equation-005}) by $e^{-i(k+1)\theta}$ and then integrating from $0$ to $2\pi$, we obtain
\begin{equation}\label{Rohi-Vasu-P3-equation-007}
c_{k} r^{k+1}G_k(r,\alpha)=\frac{1}{2\pi}\int_0^{2\pi}f_{\bar{z}}(re^{i\theta}) e^{-i(k+1)\theta}d\theta.
\end{equation}
Since $\Lambda_f(z)\leq \Lambda$, we may conclude from equations (\ref{Rohi-Vasu-P3-equation-006}) and (\ref{Rohi-Vasu-P3-equation-007}) that 
\begin{align*}\label{Rohi-Vasu-P3-equation-008}
\left(|c_kG_k(r,\alpha)|+|c_{-k}G_k(r,\alpha)|\right)r^{k+1}&\leq \frac{1}{2\pi}\int_0^{2\pi}\left(|f_z(re^{i\theta})|+|f_{\bar{z}}(re^{i\theta})|\right)d\theta\\
&=\frac{1}{2\pi}\int_0^{2\pi}\Lambda_f(re^{i\theta})d\theta\leq \Lambda.\numberthis
\end{align*}
By letting $r\to 1^{-}$ and using (\ref{Rohi-Vasu-P3-equation-003}), we obtain
$$
|c_kG_k(1,\alpha)|+|c_{-k}G_k(1,\alpha)|\leq \Lambda.
$$
This completes the proof.
\end{pf}

We obtain the following corollary from Theorem \ref{Rohi-Vasu-P3-Theorem-001}.

\begin{cor}\label{Rohi-Vasu-P3-Corollary-001}
Let $\alpha>0$ and not an even integer. Let $f\in C^2(\ID)$ be the solution of (\ref{Rohi-Vasu-P3-equation-001}) with series expansion of the form (\ref{Rohi-Vasu-P3-equation-003-a}) and $\Lambda_f(z)\leq \Lambda$ for each $z\in \ID.$ Then, for $k\in \{1,2,3,...\}$,
\begin{equation}\label{Rohi-Vasu-P3-equation-009}
|c_k|+|c_{-k}|\leq \frac{2\Lambda \Gamma(1+\frac{\alpha}{2})\Gamma(k+1+\frac{\alpha}{2})}{k! \Gamma(\alpha+1)|2k-\alpha|}.
\end{equation}
\end{cor}

\begin{pf}
Since $\alpha>0$ and $\alpha$ is not a positive even integer, by using Theorem \ref{Rohi-Vasu-P3-Theorem-001} and (\ref{Rohi-Vasu-P3-equation-003}), we obtain
the following inequality
$$
\frac{\Gamma(k+2) \Gamma(\alpha)}{\Gamma(k+1+\frac{\alpha}{2}) \Gamma(1+\frac{\alpha}{2})})\left(|c_k|+|c_{-k}|\right) \leq \frac{2\Lambda(k+1)}{\alpha|2k-\alpha|}.
$$
This implies that
$$
|c_k|+|c_{-k}|\leq \frac{2\Lambda \Gamma(1+\frac{\alpha}{2})\Gamma(k+1+\frac{\alpha}{2})}{k! \Gamma(\alpha+1)|2k-\alpha|},
$$
which completes the proof.
\end{pf}

\section{The Landau Type Theorem}
 In 2015, Chen and Vuorinen\cite{Chen-Vuorinen-2015} obtained the Landau type theorem for $\alpha$-harmoinc mappings when $\alpha\in (-1,0]$ and in  2021, Long and Wang \cite{Long-Wang-2021} obtained it when $\alpha\in (0,2)$. In this section, we derive another Landau type theorem for $\alpha$-harmonic mappings when $\alpha\in (0,2).$
\vspace{3mm}
 
We first prove the following lemma, which is vital to the proof of the main result.
\begin{lem}\label{Rohi-Vasu-P3-Lemma-001}
For $x\in [0,1)$, let
$$
\varphi(x)=\beta-\frac{2\Lambda}{2-\alpha}\left((2-\alpha)\frac{x^2}{1-x}+\frac{2a}{(1-x)^3}-2a+(2a-1)\frac{x^2}{1-x^2}\right),
$$
where $\alpha \in (0,2)$, $\beta>0$ and $\Lambda>0$ are constants and $a=\frac{\Gamma(1+\frac{\alpha}{2})}{\Gamma(1+\alpha)}$. Then $\varphi$ is strictly decreasing and there is a unique $x_0\in (0,1)$ such that $\varphi(x_0)=0.$
\end{lem}

\begin{pf}
It is easy to see that $1/2<a<1$ for $\alpha\in (0,2)$ (see \cite[Lemma 3.1]{Long-Wang-2021}). Let
$$
f_1(x)=(2-\alpha)\frac{x^2}{1-x},
$$  
$$
f_2(x)=\frac{2a}{(1-x)^3}-2a,
$$
and 
$$
f_3(x)=(2a-1)\frac{x^2}{1-x^2}.
$$
Clearly,
 $$
\varphi(x)= \beta-\frac{2\Lambda}{2-\alpha}\Big(f_1(x)+f_2(x)+f_3(x)\Big).
 $$
A simple computation shows that $f_1'(x)>0, f_2'(x)>0 $ and $f_3'(x)>0$ for $x \in (0,1)$. This implies that $\varphi(x)$ is strictly decreasing in $[0,1)$.
\vspace{2mm}

Furthermore, we can observe that
$$
\lim_{x\to 0} \varphi(x)=\beta>0 \quad \text{and}\quad \lim_{x\to 1^{-}}\varphi(x)=-\infty.
$$
Therefore, there exists a unique $x_0\in (0,1)$ such that $\varphi(x_0)=0.$ This completes the proof of this lemma.
\end{pf}

\begin{thm}\label{Rohi-Vasu-P3-Theorem-002}
For $\alpha\in (0,2)$, let $a=\Gamma(\alpha/2+1)/\Gamma(1+\alpha)$,  $f\in C^2(\ID)$ be the solution of (\ref{Rohi-Vasu-P3-equation-001}) satisfying $f(0)=0, \lambda_f(0)\geq \beta$ and $\Lambda_f(z)\leq \Lambda$ for each $z\in \ID$, where $\Lambda$, $\beta$ are positive constants. Then $f$ is univalent on a disc $\ID_{\rho_0}$, where $\rho_0$ satisfies the following equation
$$
\beta- \frac{2\Lambda}{2-\alpha}\left((2-\alpha)\frac{\rho_0^2}{1-\rho_0}+\frac{2a}{(1-\rho_0)^3}-2a
+(2a-1)\frac{\rho_0^2}{1-\rho_0^2}\right)=0. 
$$

Morever, $f(\ID_{\rho_0})$ contains a schlicht disc of radius
$$
R_0\geq \rho_0\left(\beta-\frac{2\Lambda}{2-\alpha}\left(\frac{2-\alpha}{1-\rho_0} \frac{{\rho_0}^2}{3} +\frac{2a}{(1-\rho_0)^3}-2a + \frac{2a-1}{1-{\rho_0}^2}\frac{{\rho_0}^2}{3}\right)\right).
$$
\end{thm}

\begin{pf}
By Theorem \ref{Thm-A}, we may assume that
$$
\begin{aligned}
f(z)= & \sum_{k=0}^{\infty} c_k F\left(-\frac{\alpha}{2}, k-\frac{\alpha}{2} ; k+1 ;|z|^2\right) z^k \\
& +\sum_{k=1}^{\infty} c_{-k} F\left(-\frac{\alpha}{2}, k-\frac{\alpha}{2} ; k+1 ;|z|^2\right) \bar{z}^k,\  z \in \mathbb{D},
\end{aligned}
$$
for some sequence $\left\{c_k\right\}_{k=-\infty}^{\infty}$ of complex numbers satisfying
$$
\lim _{|k| \rightarrow \infty} \sup \left|c_k\right|^{\frac{1}{k \mid}} \leq 1 .
$$
For $\alpha\in(0,2)$ and $k\in \{1,2,3,...\}$, we observe that
$$
\frac{d}{dt}F\left(-\frac{\alpha}{2}, k-\frac{\alpha}{2} ; k+1 ; t\right)=\frac{\left(-\frac{\alpha}{2}\right)\left(k-\frac{\alpha}{2}\right)}{k+1}F\left(1-\frac{\alpha}{2}, k+1-\frac{\alpha}{2} ; k+2 ; t\right)<0.
$$
This shows that $F\left(-\frac{\alpha}{2}, k-\frac{\alpha}{2} ; k+1 ; t\right)$ is decreasing on $t \in[0,1)$. Therefore, we have
\begin{equation}\label{Rohi-Vasu-P3-equation-010}
F\left(-\frac{\alpha}{2}, k-\frac{\alpha}{2} ; k+1 ; t\right) \leq F\left(-\frac{\alpha}{2}, k-\frac{\alpha}{2} ; k+1 ; 0\right)=1 .
\end{equation}
Further, for $\alpha\in(0,2)$ and $k\in\{1,2,3,...\}$, we have
\begin{align*}
\frac{d}{dt}&F\left(-\frac{\alpha}{2}+1, k-\frac{\alpha}{2}+1 ; k+2 ;t\right)\\ &=\frac{(-\frac{\alpha}{2}+1)( k-\frac{\alpha}{2}+1)}{k+2}F\left(-\frac{\alpha}{2}+2, k-\frac{\alpha}{2}+2;k+3,t \right)>0
\end{align*}
which implies that $F\left(-\frac{\alpha}{2}+1, k-\frac{\alpha}{2}+1 ; k+2 ;t\right)$ is increasing function on $[0,1)$. 
\vspace{2mm}

Note that the function $F\left(-\frac{\alpha}{2}+1, k-\frac{\alpha}{2}+1 ; k+2 ;t\right)$ is nonnegative as its parameters are nonnegative. Therefore, we have
$$
\left|F(-\frac{\alpha}{2}+1, k-\frac{\alpha}{2}+1 ; k+2 ;t)\right|\leq  \left|F(-\frac{\alpha}{2}+1, k-\frac{\alpha}{2}+1 ; k+2 ;1)\right|.
$$
From Theorem \ref{Rohi-Vasu-P3-Theorem-001}, we have
\begin{align*}\label{Rohi-Vasu-P3-equation-011}
\frac{\alpha(2k-\alpha)}{2(k+1)}&(|c_k|+|c_{-k}|)|F(-\frac{\alpha}{2}+1, k-\frac{\alpha}{2}+1 ; k+2 ;t)|\\
& \leq \frac{\alpha(2k-\alpha)}{2(k+1)}(|c_k|+|c_{-k}|)|F(-\frac{\alpha}{2}+1, k-\frac{\alpha}{2}+1 ; k+2 ;1)|\leq \Lambda.\numberthis
\end{align*}
For $k\in \{1,2,3,...\}$ and $\alpha\in (0,2)$, Long and Wang \cite{Long-Wang-2021} have proved the following inequality. 
\begin{equation}\label{Rohi-Vasu-P3-equation-012}
\left|\frac{(-\frac{\alpha}{2})_n(k-\frac{\alpha}{2})_n}{(k+1)_n}\frac{1}{n!}\right|\leq 1-\frac{\Gamma(k+1)\Gamma(1+\alpha)}{\Gamma(k+1+\frac{\alpha}{2})\Gamma(1+\frac{\alpha}{2})}.
\end{equation}
Using (\ref{Rohi-Vasu-P3-equation-012}) and Corollary \ref{Rohi-Vasu-P3-Corollary-001}, we obtain
\begin{align*}\label{Rohi-Vasu-P3-equation-013}
\left(|c_1|+|c_{-1}|\right)\left|\frac{(-\frac{\alpha}{2})_n(k-\frac{\alpha}{2})_n}{(k+1)_n}\frac{1}{n!}\right| & \leq \frac{2\Lambda}{2-\alpha}\left(\frac{
\Gamma(1+\frac{\alpha}{2})\Gamma(2+\frac{\alpha}{2})}{\Gamma(2)\Gamma(\alpha+1)}-1\right)\\
& < \frac{2\Lambda}{2-\alpha}(2a-1), \numberthis
\end{align*} 
where $a=\Gamma(\alpha/2+1)/\Gamma(1+\alpha)$ and $\alpha \in (0,2).$
 Since $f(0)=c_0=0$, a simple computation shows that
 \begin{align*}\label{Rohi-Vasu-P3-equation-014}
 f_z(z)-f_z(0) =& \sum_{k=1}^{\infty} c_k \frac{-\frac{\alpha}{2}\left(k-\frac{\alpha}{2}\right)}{k+1} F\left(-\frac{\alpha}{2}+1 ; k-\frac{\alpha}{2}+1 ; k+2 ;|z|^2\right)\bar{z} z^k\\
 &+\sum_{k=1}^{\infty} c_{-k} \frac{-\frac{\alpha}{2}\left(k-\frac{\alpha}{2}\right)}{k+1} F\left(-\frac{\alpha}{2}+1 ; k-\frac{\alpha}{2}+1 ; k+2 ;|z|^2\right)\bar{z}^{k+1}\\
 &+ \sum_{k=2}^{\infty} kc_k F\left(-\frac{\alpha}{2} ; k-\frac{\alpha}{2} ; k+1 ;|z|^2\right)z^{k-1}\\
 &+ c_1\left(F(-\frac{\alpha}{2} ; 1-\frac{\alpha}{2} ;2;|z|^2)-1\right)\numberthis
 \end{align*}
and
\begin{align*}\label{Rohi-Vasu-P3-equation-015}
 f_{\bar{z}}(z)-f_{\bar{z}}(0) =& \sum_{k=1}^{\infty} c_k \frac{-\frac{\alpha}{2}\left(k-\frac{\alpha}{2}\right)}{k+1} F\left(-\frac{\alpha}{2}+1 ; k-\frac{\alpha}{2}+1 ; k+2 ;|z|^2\right) z^{k+1}\\
 &+\sum_{k=1}^{\infty} c_{-k} \frac{-\frac{\alpha}{2}\left(k-\frac{\alpha}{2}\right)}{k+1} F\left(-\frac{\alpha}{2}+1 ; k-\frac{\alpha}{2}+1 ; k+2 ;|z|^2\right)z\bar{z}^{k}\\
 &+ \sum_{k=2}^{\infty} kc_{-k} F\left(-\frac{\alpha}{2} ; k-\frac{\alpha}{2} ; k+1 ;|z|^2\right)\bar{z}^{k-1}\\
 &+ c_{-1}\left(F(-\frac{\alpha}{2} ; 1-\frac{\alpha}{2} ;2;|z|^2)-1\right).\numberthis
\end{align*}
For $|z|=r$, in view of (\ref{Rohi-Vasu-P3-equation-014}) and (\ref{Rohi-Vasu-P3-equation-015}), we obtain
\begin{align*}
 |f_z(z)-f_z(0)|&+|f_{\bar{z}}(z)-f_{\bar{z}}(0)| \\
  \leq 2 & \sum_{k=1}^{\infty}(|c_k|+|c_{-k}|)\frac{\alpha(2k-\alpha)}{2(k+1)}\left|F\left(-\frac{\alpha}{2}+1, k-\frac{\alpha}{2}+1 ; k+2 ;r^2\right)\right|r^{k+1}\\
 &+ \sum_{k=2}^{\infty}k(|c_k|+|c_{-k}|)\left|F\left(-\frac{\alpha}{2} ; k-\frac{\alpha}{2} ;k+1;r^2\right)\right|r^{k-1}\\
 &+ (|c_1|+|c_{-1}|)\left|F\left(-\frac{\alpha}{2} ; 1-\frac{\alpha}{2} ;2;r^2\right)-1\right|. 
\end{align*}
Applying (\ref{Rohi-Vasu-P3-equation-010}), (\ref{Rohi-Vasu-P3-equation-011}), (\ref{Rohi-Vasu-P3-equation-013}) and Corollary \ref{Rohi-Vasu-P3-Corollary-001}, we obtain
\begin{align*}\label{Rohi-Vasu-P3-equation-016}
|f_z(z)&-f_z(0)|+|f_{\bar{z}}(z)-f_{\bar{z}}(0)|\\
& \leq  2\Lambda \sum_{k=1}^{\infty} r^{k+1}+ \sum_{k=2}^{\infty}k(|c_k|+|c_{-k}|)r^{k-1}+ \left(|c_1|+|c_{-1}|\right)\sum_{n=1}^{\infty}\left|\frac{(-\frac{\alpha}{2})_n(k-\frac{\alpha}{2})_n}{(k+1)_n}\frac{r^{2n}}{n!}\right|\\
&\leq 2\Lambda \sum_{k=1}^{\infty} r^{k+1}+\sum_{k=2}^{\infty}\frac{2\Lambda a \Gamma(k+1+\frac{\alpha}{2})}{\Gamma(k+1) (2k-\alpha)}r^{k-1}+\frac{2\Lambda}{2-\alpha}(2a-1)\sum_{n=1}^{\infty}r^{2n}\\
& < 2\Lambda \sum_{k=1}^{\infty} r^{k+1}+\frac{2\Lambda}{2-\alpha}a\sum_{k=2}^{\infty}k(k+1)r^{k-1}+\frac{2\Lambda}{2-\alpha}(2a-1)\sum_{n=1}^{\infty}r^{2n}\\
&=2\Lambda\frac{r^2}{1-r}+\frac{2\Lambda}{2-\alpha} \left(\frac{2a}{(1-r)^3}-2a\right)+\frac{2\Lambda}{2-\alpha}(2a-1)\frac{r^2}{1-r^2}\\
&=\frac{2\Lambda}{2-\alpha}\left((2-\alpha)\frac{r^2}{1-r}+\frac{2a}{(1-r)^3}-2a+(2a-1)\frac{r^2}{1-r^2}\right). \numberthis
\end{align*}
To prove the univalence of $f$ in $\ID_{\rho_0}$, we choose two distinct points $z_1,z_2 \in \ID_{\rho_0}$ and let $[z_1,z_2]$ denote the line segment from $z_1$ to $z_2$ with ends points $z_1$ and $z_2$, where $\rho_0$ satisfies the following equation
$$
\beta- \frac{2\Lambda}{2-\alpha}\left((2-\alpha)\frac{\rho_0^2}{1-\rho_0}+\frac{2a}{(1-\rho_0)^3}-2a
+(2a-1)\frac{\rho_0^2}{1-\rho_0^2}\right)=0.
$$  
By using (\ref{Rohi-Vasu-P3-equation-016}) and Lemma \ref{Rohi-Vasu-P3-Lemma-001}, we have
\begin{align*}
\left|f\left(z_2\right)-f\left(z_1\right)\right|= & \left|\int_{\left[z_1, z_2\right]} f_z(z) d z+f_{\bar{z}}(z) d \bar{z}\right| \\
= & \left|\int_{\left[z_1, z_2\right]} f_z(0) d z+f_{\bar{z}}(0) d \bar{z}\right| \\
& -\left|\int_{\left[z_1, z_2\right]}\left(f_z(z)-f_z(0)\right) d z+\left(f_{\bar{z}}(z)-f_{\bar{z}}(0)\right) d \bar{z}\right| \\
\geq & \lambda_f(0)\left|z_2-z_1\right| \\
& -\int_{\left[z_1, z_2\right]}\left(\left|f_z(z)-f_z(0)\right|+\left|f_{\bar{z}}(z)-f_{\bar{z}}(0)\right|\right)|d z|\\
> & |z_1-z_2| \left\{\beta- \frac{2\Lambda}{2-\alpha}\left((2-\alpha)\frac{\rho_0^2}{1-\rho_0}+\frac{2a}{(1-\rho_0)^3}-2a \right. \right. \\
&+ \left. \left. (2a-1)\frac{\rho_0^2}{1-\rho_0^2}\right)\right\}\\
&=0. 
\end{align*}
Therefore, $f(z_1)\neq f(z_2)$. This shows that $f$ is univalent in the ball $\D_{\rho_0}.$
\vspace{3mm}

To prove the second part of the theorem, let $ \zeta = \rho_0 e^{i\theta}\in \partial \D_{\rho_0} $, then we have 
\begin{align*}
|f(\zeta)-f(0)| = & \left|\int_{[0, \zeta]} f_z(z) d z+f_{\bar{z}}(z) d \bar{z}\right| \\
 = &\left|\int_{[0, \zeta]} f_z(0) d z+f_{\bar{z}}(0) d \bar{z}\right| \\
&-\left|\int_{[0, \zeta]}\left(f_z(z)-f_z(0)\right) d z+\left(f_{\bar{z}}(z)-f_{\bar{z}}(0)\right) d \bar{z}\right| \\
\geq & \lambda_f(0) |\zeta|-\int_{[0, \zeta]}\big(\left|f_z(z)-f_z(0)\right|+\left|f_{\bar{z}}(z)-f_{\bar{z}}(0)\right|\big)|d z|\\
\geq & \beta \rho_0 - \frac{2\Lambda}{2-\alpha}\int_{0}^{\rho_0} \left((2-\alpha)\frac{r^2}{1-r}+\frac{2a}{(1-r)^3}-2a \right.\\ 
& + \left. (2a-1)\frac{r^2}{1-r^2}\right)dr\\
\geq & \beta\rho_0-\frac{2\Lambda}{2-\alpha}\left(\frac{2-\alpha}{1-\rho_0}\int_{0}^{\rho_0} r^2 dr + \left(\frac{2a}{(1-\rho_0)^3}-2a\right)\int_{0}^{\rho_0}dr \right.\\
&+ \left. \frac{2a-1}{1-{\rho_0}^2} \int_{0}^{\rho_0} r^2 dr\right)\\
= & \beta\rho_0-\frac{2\Lambda}{2-\alpha}\left(\frac{2-\alpha}{1-\rho_0} \frac{{\rho_0}^3}{3} + \left(\frac{2a}{(1-\rho_0)^3}-2a\right)\rho_0 + \frac{2a-1}{1-{\rho_0}^2}\frac{{\rho_0}^3}{3}\right)\\
=& \rho_0\left(\beta-\frac{2\Lambda}{2-\alpha}\left(\frac{2-\alpha}{1-\rho_0} \frac{{\rho_0}^2}{3} +\frac{2a}{(1-\rho_0)^3}-2a + \frac{2a-1}{1-{\rho_0}^2}\frac{{\rho_0}^2}{3}\right)\right).
\end{align*}
Hence $f(\ID_{\rho_0})$ contains a Schlicht disc of radius $R_0$ with
$$
R_0\geq \rho_0\left(\beta-\frac{2\Lambda}{2-\alpha}\left(\frac{2-\alpha}{1-\rho_0} \frac{{\rho_0}^2}{3} +\frac{2a}{(1-\rho_0)^3}-2a + \frac{2a-1}{1-{\rho_0}^2}\frac{{\rho_0}^2}{3}\right)\right).
$$
This completes the proof.
\end{pf}

\begin{cor}\label{Rohi-Vasu-P3-Corollary-002}
For $\alpha\in (0,2)$, let $a=\Gamma(\alpha/2+1)/\Gamma(1+\alpha)$,  $f\in C^2(\ID)$ be the solution of (\ref{Rohi-Vasu-P3-equation-001}) satisfying $f(0)=|J_f(0)|-\beta=0$ and $\Lambda_f(z)\leq \Lambda$ for each $z\in \ID$, where $\Lambda$, $\beta$ are positive constants and $J_f$ is the Jaccobian of $f$. Then $f$ is univalent on a disc $\ID_{\rho_1}$, where $\rho_1$ satisfies the following equation
$$
\frac{\beta}{\Lambda}- \frac{2\Lambda}{2-\alpha}\left((2-\alpha)\frac{\rho_1^2}{1-\rho_1}+\frac{2a}{(1-\rho_1)^3}-2a
+(2a-1)\frac{\rho_1^2}{1-\rho_1^2}\right)=0.
$$ 
Morever, $f(\ID_{\rho_1})$ contains a schlicht disc of radius
$$
R_1\geq \rho_1\left(\frac{\beta}{\Lambda}-\frac{2\Lambda}{2-\alpha}\left(\frac{2-\alpha}{1-\rho_1} \frac{{\rho_1}^2}{3} +\frac{2a}{(1-\rho_1)^3}-2a + \frac{2a-1}{1-{\rho_1}^2}\frac{{\rho_1}^2}{3}\right)\right).
$$
\end{cor}
\begin{pf}
Since, $|J_f(0)|=\beta$, we have
$$
\beta=|J_f(0)|=|\operatorname{det} D_f(0)|=\Lambda_f(0)\lambda_f(0)\leq \Lambda\lambda_f(0),
$$
which gives
$$
\lambda_f(0)=|| f_z(0)|-| f_{\bar{z}}(0)||\geq \frac{\beta}{\Lambda}.
$$
Let $z_1, z_2 \in \ID_{\rho_1}$, then by using the proof of Theorem \ref{Rohi-Vasu-P3-Theorem-002} and Lemma \ref{Rohi-Vasu-P3-Lemma-001}, we have
\begin{align*}
|f(z_1)-f(z_2)| > |z_1-z_2| &\left\{\frac{\beta}{\Lambda}- \frac{2\Lambda}{2-\alpha}\left((2-\alpha)\frac{\rho_1^2}{1-\rho_1}+\frac{2a}{(1-\rho_1)^3}-2a \right. \right. \\
&+ \left. \left. (2a-1)\frac{\rho_1^2}{1-\rho_1^2}\right)\right\}=0.
\end{align*}
Therefore, $f(z_1)\neq f(z_2)$. This shows that $f$ is univalent in the ball $\D_{\rho_1}.$
\vspace{2mm}

To prove the second part of the corollary, let $\zeta=\rho_1 e^{i\theta}\in \partial \ID_{\rho_1}$. By using the proof of Theorem \ref{Rohi-Vasu-P3-Theorem-002}, we obtain
\begin{align*}
|f(\zeta)|\geq \rho_1\left(\frac{\beta}{\Lambda}-\frac{2\Lambda}{2-\alpha}\left(\frac{2-\alpha}{1-\rho_1} \frac{{\rho_1}^2}{3} +\frac{2a}{(1-\rho_1)^3}-2a + \frac{2a-1}{1-{\rho_1}^2}\frac{{\rho_1}^2}{3}\right)\right).
\end{align*}
Hence $f(\ID_{\rho_1})$ contains a Schlicht disc of radius $R_1$ with
$$
R_1\geq  \rho_1\left(\frac{\beta}{\Lambda}-\frac{2\Lambda}{2-\alpha}\left(\frac{2-\alpha}{1-\rho_1} \frac{{\rho_1}^2}{3} +\frac{2a}{(1-\rho_1)^3}-2a + \frac{2a-1}{1-{\rho_1}^2}\frac{{\rho_1}^2}{3}\right)\right).
$$
This completes the proof.

\end{pf}

\noindent\textbf{Acknowledgement:}  The second named author thank CSIR for their support. 

\end{document}